\documentclass[11pt]{amsart}
\usepackage[a4paper,hmargin={2.25cm,2.25cm},vmargin={2.25cm,2.25cm}]{geometry}
\usepackage[english]{babel} 
\usepackage[utf8]{inputenc} 
\usepackage[T1]{fontenc}
\usepackage{lmodern} 
\usepackage{thmtools} 
\usepackage{mathtools}
\usepackage{mleftright}
\usepackage{amssymb} 
\usepackage{hyperref}
\usepackage{mleftright}
\usepackage[shortlabels]{enumitem}
\usepackage{tikz-cd}
\tikzset{%
  symbol/.style={%
    draw=none, every to/.append style={%
      edge node={node [sloped, allow upside down,auto=false]{$#1$}}} } }

\usepackage{floatrow}
\usepackage{graphicx}
\usepackage{caption}
\usepackage{stmaryrd}
\usepackage{tikz}
\usetikzlibrary{arrows}
\usepackage{todonotes}
\usepackage{mleftright}
\usepackage{cleveref} 

\declaretheorem[
name         = Theorem,
refname      = {Theorem,Theorems},
Refname      = {Theorem,Theorems},
numberwithin = section,
]{theorem}

\declaretheorem[
name         = Proposition,
refname      = {Proposition,Propositions},
Refname      = {Proposition,Propositions},
sibling      = theorem,
]{proposition}

\declaretheorem[
name         = Lemma,
refname      = {Lemma,Lemmas},
Refname      = {Lemma,Lemmas},
sibling      = theorem,
]{lemma}

\declaretheorem[
name         = Corollary,
refname      = {Corollary,Corollaries},
Refname      = {Corollary,Corollaries},
sibling      = theorem,
]{corollary}

\theoremstyle{definition}

\declaretheorem[
name        = Definition,
refname     = {Definition,Definitions},
Refname     = {Definition,Definitions},
sibling     = theorem,
]{definition}

\declaretheorem[
name        = Example,
refname     = {Example,Examples},
Refname     = {Example,Examples},
sibling     = theorem,
]{example}

\declaretheorem[
name        = Remark,
refname     = {Remark,Remarks},
Refname     = {Remark,Remarks},
sibling     = theorem,
]{remark}

\newlist{tfae}{enumerate}{1}%
\setlist[tfae,1]{label=(\roman*)}%

\makeatletter
\def\namedlabel#1#2{\begingroup
    #2%
    \def\@currentlabel{#2}%
    \phantomsection\label{#1}\endgroup
}

\newcommand{\cat}[1]{\mathsf{#1}}

\newcommand{\op}{\mathrm{op}}
\newcommand{\MetCHsep}{\cat{MetCH_{sep}}}
\newcommand{\MetCHsepop}{\cat{MetCH}_\cat{sep}^\op}
\newcommand{\MetCHsepsym}{\cat{MetCH_{sep,sym}}}
\newcommand{\MetCHsepsymop}{\cat{MetCH_{sep,sym}^\op}}
\newcommand{\CH}{\cat{CH}}

\newcommand{\Set}{\cat{Set}}

\newcommand{\N}{\mathbb{N}}
\newcommand{\R}{\mathbb{R}}

\title[The dual of separated metric compact Hausdorff spaces is a variety]{The dual of separated metric compact Hausdorff spaces \\is an infinitary variety of algebras}

\author{Marco Abbadini}

\address[Marco Abbadini]{Research Institute in Mathematics and Physics, Université catholique de Louvain,
Chem.\ du Cyclotron 2, 1348 Ottignies-Louvain-la-Neuve, Belgium}
\urladdr{\url{https://marcoabbadini-uni.github.io/}}

\author{Dirk Hofmann}

\address[Dirk Hofmann]{Center for Research and Development in Mathematics and Applications, Department of Mathematics, University of Aveiro, Portugal}
\urladdr{\url{https://sweet.ua.pt/dirk/}}

\thanks{The first author was funded by an FSR Incoming Postdoctoral Fellowship. The second author acknowledges the partial support by CIDMA under the Portuguese Foundation for Science and Technology (FCT, \url{https://ror.org/00snfqn58}) Multi-Annual Financing Program for R\&D Units, grants \url{UID/4106/2025} and \url{UID/PRR/4106/2025}.}



\begin{document}

\begin{abstract}
  A \emph{metric compact Hausdorff space} is a Lawvere metric space equipped with a \emph{compatible} compact Hausdorff topology (which does not need to be the induced topology).
  These spaces maintain many important features of compact metric spaces, but the resulting category is much better behaved.
	We prove that the category of separated metric compact Hausdorff spaces is dually equivalent to a variety of infinitary algebras.
\end{abstract}

\maketitle

\tableofcontents

\section{Introduction}
\label{sec:introduction}

Stone duality \cite{Sto36} identifies the category of Stone spaces with the dual of the category of Boolean algebras. Since Boolean algebras form a variety, this gives an algebraic description of these spaces in terms of operations and equations. Many subsequent dualities follow this pattern. In particular, Duskin \cite[5.15.3]{Dus69} and Negrepontis \cite{Neg71} showed that compact Hausdorff spaces are dually equivalent to a variety of infinitary algebras; in categorical terms, this means that the opposite of the category of compact Hausdorff spaces is monadic over $\Set$ (see also \cite[Ch.~9, Cor.~1.11]{BarrWells2005} for the missing details in \cite{Dus69}).

Priestley duality \cite{Pri70,Pri72} provides an ordered counterpart: Priestley spaces are dually equivalent to bounded distributive lattices, which again form a variety. More recently, this algebraic description was extended to Nachbin's compact ordered spaces \cite{Nac65}. Following the quasivariety result of \cite{HNN18}, it was proved in \cite{Abb19,AR20} that the category of Nachbin's compact ordered spaces is also dually equivalent to a variety of infinitary algebras.

Lawvere's interpretation of orders and metrics as enriched categories \cite{Law73} suggests a further step: replacing order by metric. Here a metric is a map $d\colon X\times X\to[0,\infty]$ satisfying $d(x,x)=0$ and the triangle inequality, without any assumption of symmetry. It is \emph{separated} if $d(x,y)=d(y,x)=0$ implies $x=y$. The resulting topological structures are \emph{separated metric compact Hausdorff spaces}: compact Hausdorff spaces equipped with a lower semicontinuous separated metric (whose induced topology need not be the given topology). These include Nachbin's compact ordered spaces, viewed as spaces with $\{0,\infty\}$-valued metrics, as well as classical compact metric spaces. We denote by $\MetCHsep$ the category of separated metric compact Hausdorff spaces and non-expansive continuous maps between them.

The preceding dualities lead naturally to the question of whether $\MetCHsep$ is likewise dually equivalent to a variety. In \cite{AH25}, we proved that $\MetCHsep$ is complete, cocomplete and Barr-coexact, but left the question open \cite[Question 5.18]{AH25}. By a standard categorical characterisation of varieties \cite[Thm.~4.4.5]{Bor94a}, these properties reduce the remaining task to proving the existence of a regular injective regular cogenerator of $\MetCHsep$. In this paper we exhibit such an object, thereby proving that $\MetCHsepop$ is equivalent to a variety of infinitary algebras.

The main tool is a Tietze--Urysohn type extension theorem, essentially due to Matou\v{s}kov\'a \cite[Thm.~2.4]{Mat00}, of which we were unaware when writing \cite{AH25}. Adapting her argument (and closing a minor gap therein) to asymmetric metrics and possibly infinite distances yields the extension property needed here: every non-expansive continuous map from a closed subspace into a bounded closed real interval extends to the whole space. The interval carries its usual topology and the asymmetric metric $d(u,v)=\max\{v-u,0\}$. Thus these intervals are regular injective. The same extension property also shows that
\[
  \prod_{n\in\N}[0,n]
\]
is a regular cogenerator. Since it is a product of regular injective objects, it is also regular injective. Thus, $\prod_{n\in\N}[0,n]$ supplies the required regular injective regular cogenerator object of $\MetCHsep$.

We use the extension theorem also to answer the question after \cite[Cor.~2.22]{HN23}: we show that every separated metric compact Hausdorff space is Priestley (in the metric sense, see \cite[Def.~2.18]{HN23}).

The symmetric objects of $\MetCHsep$ are the compact topometric spaces of \cite{BenYaacov2008,BenYaacov2013}, which include the type spaces of continuous first-order logic \cite{BenYaacovUsvyatsov2010}. Ben Yaacov \cite[Thms.~1.6 and~1.9]{BenYaacov2013} characterised topometric normality in terms of a Lipschitz version of Urysohn's lemma and a Tietze extension theorem allowing an arbitrarily small increase in the Lipschitz bound.

These spaces are also precisely the compact extended metric-topological spaces \cite[Rem.~3.2]{PST26}. Ambrosio, Erbar and Savar\'e introduced extended metric-topological spaces for optimal transport and infinite-dimensional metric analysis, with abstract Wiener spaces as a motivating example \cite{AES16}. For these compact spaces, Pasqualetto and Taipalus \cite[Thm.~3.1 and Rem.~3.2(i)]{PT26} extend Matou\v{s}kov\'a's theorem to possibly infinite distances by truncating the metric.

\section{Preliminaries}
\label{sec:preliminaries}

We recall from the introduction the following definitions.

\begin{definition}[(Separated) metric, non-expansive map] \hfill
	\begin{enumerate}

		\item
		By a \emph{metric} on a set $X$ we mean a map \(d \colon X\times X\to [0,\infty]\) such that, for all \(x,y,z\in X\),
  \begin{displaymath}
    d(x,x)=0
    \quad\text{and}\quad
    d(x,z)\leq d(x,y)+d(y,z).
  \end{displaymath}
  We say that \(d\) is \emph{separated} if, for all \(x,y\in X\), \(d(x,y)=0=d(y,x)\) implies \(x=y\).

  	\item
	A map $f \colon X\to X'$ between metric spaces $(X,d)$ and $(X',d')$ is called \emph{non-expansive} if $d'(f(x),f(y))\leq d(x,y)$ for all $x,y\in X$.

	\end{enumerate}
\end{definition}

  The following is the category of interest in this paper. Recall that the \emph{upper topology} on \([0,\infty]\) is generated by the sets \(]u,\infty]\).
\begin{definition}[\cite{HR18}]
  \label{d:def:1}
  A \emph{(separated) metric compact Hausdorff space} \(X\) is a compact Hausdorff space \(X\) together with a (separated) metric  $d \colon X \times X \to [0,\infty]$ that is lower semicontinuous, i.e., continuous with respect to the upper topology of $[0,\infty]$.

	We denote the category of separated metric compact Hausdorff spaces and non-expansive continuous maps by \(\MetCHsep\).
\end{definition}

\begin{example}\label{d:ex:1}
  Let $a,b \in \R$ with $a \leq b$. The interval \([a,b]\) becomes a separated metric space with the Euclidean topology and the metric \(d\) defined by $d(x,y) = \max(y-x, 0)$.
  Its dual, \([a,b]^{\circ}\), has distance $d^{\circ}(x,y) = \max(x-y, 0)$ and is isomorphic to \([a,b]\) via the involution \(t\mapsto a+b-t\).
\end{example}

In this paper we wish to disclose the algebraic nature of the dual category of \(\MetCHsep\).
We will make use of the following characterisation of varieties, which can be found, for instance, in \cite[Thm.~4.4.5]{Bor94a}.

\begin{theorem}
  \label{d:thm:1}
  A category is equivalent to a (possibly infinitary) variety (i.e., is monadic over $\Set$) if and only if it is complete, cocomplete and Barr-exact, and has a regular projective regular generator.
\end{theorem}

We collect some properties that we already know of $\MetCHsep$ and that will be useful to show its dual to be equivalent to a variety.

An \emph{embedding} $\iota \colon X \hookrightarrow Y$ between metric compact Hausdorff spaces is an injective continuous map that preserves distances; since the spaces are compact Hausdorff, $\iota$ is automatically a homeomorphism onto its image.

\begin{theorem}
  \label{p:regmono=embeddings}
  In $\MetCHsep$, the surjective morphisms are precisely the epimorphisms and the embeddings are precisely the regular monomorphisms.
\end{theorem}
\begin{proof}
  See \cite[Props.~4.5 and~4.6]{AH25}.
\end{proof}

\begin{remark}
  \label{r:reg-cogen-criterion}
  Hence, a separated metric compact Hausdorff space \(G\) with distance \(\delta\) is a regular cogenerator of \(\MetCHsep\) if and only if, for each separated metric compact Hausdorff space \(X\) with metric \(d\) and all \(x,y\in X\),
  \begin{displaymath}
    d(x,y)=\sup\{\delta(f(x),f(y))\mid f \colon X\to G\text{ in }\MetCHsep\}.
  \end{displaymath}
\end{remark}

We say that a category is \emph{Barr-coexact} if its opposite is Barr-exact.

\begin{theorem}
  $\MetCHsep$ is complete, cocomplete and Barr-coexact.
\end{theorem}
\begin{proof}
  See \cite[Thms.~2.5 and~5.16]{AH25}.
\end{proof}

\section{Algebraicity of \texorpdfstring{$\MetCHsepop$}{MetCH-sep (opposite)}}

For our main result we will use a topological-metric version of the Tietze extension theorem.
This is essentially due to Matou\v{s}kov\'a \cite[Thm.~2.4]{Mat00}, who provided this result for separated symmetric finite-valued metrics.
We will observe that her proof goes through without assuming separation and symmetry, and the case of possibly infinite-valued metrics follows immediately.
We also fix a slight gap in her argumentation.

Given a metric space $X$ with metric $d$, $A \subseteq X$ and $x \in X$, we write
\begin{displaymath}
    d(A, x) \coloneqq \inf\{d(a, x) \mid a \in A\}
  \qquad\text{and}\qquad
    d(x, A) \coloneqq \inf\{d(x, a) \mid a \in A\};
\end{displaymath}
here the infimum is computed in $[0, \infty]$. In the following definition we use the nomenclature from \cite{BenYaacov2008,BenYaacov2013}; the condition appeared already without a name in \cite{Mat00}.

\begin{definition}[Closed metric neighbourhoods]
  Let $X$ be a topological space equipped with a metric $d$. We say that \emph{$X$ admits closed metric neighbourhoods} whenever, for every closed $A\subseteq X$, the maps \(d(A,-)\colon X\to [0,\infty]\) and \(d(-,A)\colon X\to [0,\infty]\) are lower semicontinuous, i.e., for every $r>0$, both $\{x\in X\mid d(A,x)\leq r\}$ and $\{x\in X\mid d(x,A)\leq r\}$ are closed.
\end{definition}

The observation immediately preceding \cite[Lemma~2.1]{Mat00} has a small gap. It asserts that a certain metric $d'$ arising from a given topological-metric space $(X,d)$ and a non-expansive continuous map $f \colon X \to \R$ inherits the property of having closed metric neighbourhoods; however, this can fail if $f$ is unbounded.%
\footnote{Let $X$ be the topological disjoint union of a convergent sequence $C=\{x_\infty,x_1,x_2,\ldots\}$ with $x_n\to x_\infty$ and a discrete space $A=\{a_1,a_2,\ldots\}$. Give $X$ the shortest-path metric $d$ with unit edges between distinct points of $C$ and an edge of length $n$ from $a_n$ to $x_n$. For a closed $B\subseteq X$, the centres in $B\cap C$ contribute either $B\cap C$ or $C$ to the intersection of a ball about $B$ with $C$. At a fixed finite radius only finitely many centres in $B\cap A$ can contribute, each giving a singleton or all of $C$. Since the distances are integers and $A$ is a discrete summand, balls about closed sets are closed. Define $f=0$ on $C$ and $f(a_n)=n$. Then $f$ is continuous and $1$-Lipschitz. For $d_f(u,v)=d(u,v)-f(u)+f(v)$ we have $d_f(A,x_n)=0$ and $d_f(A,x_\infty)=1$, so the directed ball of radius $1/2$ about $A$ is not closed.}
Moreover, her proof does not work even with the corrected statement. (The passage that does not work is ``in both cases $x \in U_2$'' in the case $\beta \neq 0$.)
We present in the following lemma the corrected statement and proof.
The main results of \cite{Mat00} involve bounded maps, so the following version suffices.

\begin{lemma}\label{d:lem:1}
  Let $d$ be a finite-valued metric on a topological space $X$. Assume that $(X,d)$ admits closed metric neighbourhoods. Let $[a,b] \subseteq \R$ be a closed interval, and let $f\colon X\to[a,b]$ be a continuous map such that $f(x)-f(y)\leq d(x,y)$ for all $x,y\in X$. Then
  \[
    d_f(x,y)\coloneqq d(x,y)-f(x)+f(y)
  \]
  is a metric, and $(X,d_f)$ admits closed metric neighbourhoods.
\end{lemma}

\begin{proof}
  The inequality for $f$ makes $d_f$ nonnegative.
  Clearly, $d_f(x,x) = 0$ for all $x \in X$.
  The triangle inequality for $d$ gives the triangle inequality for $d_f$:
  \[
    d_f(x,z) = d(x,z) - f(x) + f(z)
    \leq d(x,y) + d(y,z) - f(x) + f(z) = d_f(x,y) + d_f(y,z).
  \]

  Let $A\subseteq X$ be closed. If $A$ is empty, both distance maps $d_f(A,-)$ and $d_f(-,A)$ are constant with value $\infty$. Otherwise, set
  \[
    W(x)\coloneqq\inf_{u\in A}\bigl(d(u,x)-f(u)\bigr).
  \]
  For each $n\in\N$, cover $[a,b]$ by finitely many closed intervals $[\ell_j,u_j]$ of length at most $2^{-n}$ and set $A_j=A\cap f^{-1}([\ell_j,u_j])$, discarding the empty sets. Each $A_j$ is closed, so
  \[
    f_n(x)\coloneqq\min_j\bigl(d(A_j,x)-u_j\bigr)
  \]
  is lower semicontinuous as a finite minimum of lower semicontinuous maps.
  For $u\in A_j$ we have $-u_j\leq -f(u)\leq -u_j+2^{-n}$.
  Taking the infimum over each $A_j$ and then the minimum over $j$ gives
  \begin{equation}\label{d:eq:1}
    f_n(x)\leq W(x)\leq f_n(x)+2^{-n}.
  \end{equation}
  By \eqref{d:eq:1}, $W=\sup_{n\in\N}f_n$, so $W$ is lower semicontinuous.
  Since $f$ is continuous, $d_f(A,-)=W+f$ is lower semicontinuous.
  Semicontinuity of $d_f(-,A)$ can be obtained by applying the same argument to the metric $d^{\circ}(x,y) \coloneqq d(y,x)$ and the function $-f$.
\end{proof}

Using now the lemma above, the proofs of \cite[Thms.~2.2 and 2.4]{Mat00} go through, giving us the following:

\begin{theorem}[A topometric Tietze extension theorem {\cite[Thm.~2.4]{Mat00}}]
  \label{t:TIetze-Urysohn}
  Let $d$ be a (possibly non-separated, possibly non-symmetric, possibly infinite-valued) metric on a normal topological space $X$ such that \(X\) admits closed metric neighbourhoods.
  Let $a,b \in \R$ be such that $a \leq b$, and let $F \subseteq X$ be closed. Every non-expansive continuous map $F \to [a,b]$  (where $[a,b]$ is equipped with the asymmetric metric, see \cref{d:ex:1}) has a non-expansive continuous extension $X \to [a,b]$.
\end{theorem}

\begin{proof}
  Matou\v{s}kov\'a proved this result under the additional hypotheses that the metric \(d\) is separated, symmetric and finite-valued. Her proof goes through (up to the small correction we made in \cref{d:lem:1}) also without separation and symmetry, essentially because her technical Lemma~2.1 does not assume these conditions.
  Let us note that she considers  the distance \(d(x,y)=\max\{0,x-y\}\) in \([a,b]\); i.e., she considers \([a,b]^{\circ}\) in our notation; however, this difference is irrelevant since \([a,b]\cong[a,b]^{\circ}\).

  To handle possibly infinite distances, we truncate the metric as in the symmetric case treated by Pasqualetto and Taipalus \cite[proof of Thm.~3.1]{PT26}.

  Given $M \in [0, \infty[$ and a metric space $Y$ with metric $d$, we write $\lfloor Y \rfloor_M$ for the metric space with the same underlying set and metric $d_M$ defined by
  \[
    d_M(x,y) = \min\{d(x,y), M\}.
  \]

  Assume now that \(X\) satisfies the conditions of \cref{t:TIetze-Urysohn}. Then, for every nonempty closed set $A$ we have
  \begin{displaymath}
    \{x \in X \mid d_M(A,x) \leq \varepsilon\}
    =\{x \in X \mid d(A,x) \leq \varepsilon\}
  \end{displaymath}
  for \(\varepsilon<M\) and
  \begin{displaymath}
    \{x \in X \mid d_M(A,x) \leq \varepsilon\}=X
  \end{displaymath}
  for \(\varepsilon\geq M\) (and similarly for $d_M(x,A)$).
  For $A=\varnothing$, both distance maps are constantly $\infty$.
 Thus also $\lfloor X \rfloor_M$ satisfies the conditions of \cref{t:TIetze-Urysohn}.

  Let \(F\) be a closed subspace of \(X\) and let $f \colon F \to [a,b]$ be a non-expansive continuous map. Set $M \coloneqq b - a$.
  Then $f \colon \lfloor F \rfloor_M \to [a,b]$ is still non-expansive, and of course continuous. By the finite-valued case of \cref{t:TIetze-Urysohn}, $f \colon \lfloor F \rfloor_M \to [a,b]$ can be extended to a non-expansive continuous map $\tilde{f} \colon \lfloor X \rfloor_M \to [a,b]$.
  The composite
  \[
    X \longrightarrow \lfloor X\rfloor_M \xrightarrow{\;\,\tilde{f}\;\,} [a,b]
  \]
  is an extension of $f \colon F \to [a,b]$ with the desired properties.
\end{proof}

\begin{corollary}[{$[a,b]$ is regular injective}]
  \label{t:reg-inj}
  Let $a,b \in \R$ with $a \leq b$.
  For every metric compact Hausdorff space \(X\) and every closed $F \subseteq X$, every non-expansive continuous map $F \to [a,b]$  (where $[a,b]$ is equipped with the asymmetric metric, see \cref{d:ex:1}) has a non-expansive continuous extension $X \to [a,b]$.
  In particular, $[a,b]$ (with the asymmetric metric) is a regular injective object in $\MetCHsep$.
\end{corollary}

\begin{proof}
  We show that every metric compact Hausdorff space satisfies the conditions of \cref{t:TIetze-Urysohn} (cf.\ also \cite[Lem.~1.9]{BenYaacov2008}).
  Clearly, $X$ is a normal topological space.
  Let $A$ be a closed subset of $X$, let $\varepsilon >0$, and let us prove that $\{x \in X \mid d(A, x) \leq \varepsilon\}$ is closed.
  We note that, for all $x \in X$,
  \begin{equation} \label{eq:distance-less-than-eps}
  d(A,x) \leq \varepsilon \iff \exists a \in A : d(a,x) \leq \varepsilon.
  \end{equation}
  The right-to-left implication is obvious, while the left-to-right one follows from the following fact: the image $B$ of $A$ under the semicontinuous map $d(-,x) \colon X \to [0, \infty]$ is compact in the space $[0, \infty]$ with the upper topology (because $A$ is compact).
 Since $\inf B = d(A,x) \leq \varepsilon $ and $B$ is compact in the space $[0, \infty]$ with the upper topology, $B$ has a minimum $m$, and $m \leq \varepsilon$.
 Thus, there is $a \in A$ such that $d(a,x) = m \leq \varepsilon $.
 This shows \eqref{eq:distance-less-than-eps}.
 Let $\pi_2 \colon X \times X \to X$ denote the projection onto the second coordinate.
 Then, by \eqref{eq:distance-less-than-eps}, $\{ x \in X \mid d(A, x) \leq \varepsilon\} = \pi_2[(A \times X) \cap d^{-1}[[0,\varepsilon]]]$, which is closed.
 A similar argument shows that $\{ x \in X \mid d(x,A) \leq \varepsilon\}$ is closed.
Since embeddings coincide with regular monomorphisms (\cref{p:regmono=embeddings}), $[a,b]$ is regular injective.
\end{proof}

\begin{remark}
  In \cref{t:reg-inj} we cannot consider an unbounded interval such as \([0,\infty]\).
  Here \([0,\infty]\) becomes a separated metric compact Hausdorff space with the metric \(d\) defined by $d(u,v) = \max(v - u, 0)$ (see \cite{Law73}, and note that \(\infty-\infty=0\)) and the \emph{Lawson topology} \cite{GHK+03}, which is generated by the subbasic open subsets $[0,u[$ and $]u,\infty]$, with $u\in [0,\infty]$.
  
  We show that \([0,\infty]\) is not regular injective.
  For example, consider the subspace \(\{0,\infty\}\) of \([0,\infty]^{\circ}\) and the non-expansive continuous map $h\colon \{0,\infty\}\to [0,\infty]$ with $h(0)=\infty$ and $h(\infty)=0$.
  If $\bar{h}\colon [0,\infty]^{\circ}\to [0,\infty]$ were a non-expansive extension of \(h\), then for every finite $t\in [0,\infty]$ we would have
  \[
    d_{[0,\infty]}(\bar{h}(t),\infty)
      \leq d_{[0,\infty]^{\circ}}(t,0)=t.
  \]
  Since $d_{[0,\infty]}(u,\infty)=\infty$ for finite $u$, this forces $\bar{h}(t)=\infty$ for every finite $t$. Such an extension cannot be continuous at $\infty$, where it takes the value $0$. Thus, $[0,\infty]$ is not regular injective.
\end{remark}

\begin{theorem} \label{t:reg-cogen}
  The product
  \[
    \prod_{n \in \N} [0,n]
  \]
  (where $[0,n]$ carries the asymmetric metric, see \cref{d:ex:1}) is a regular cogenerator in $\MetCHsep$.
\end{theorem}

      \begin{proof}
  Let \(X\in\MetCHsep\) and \(x,y\in X\). For each \(n\in\N\), the map \(g_n\colon\{x,y\}\to[0,n]\) given by \(g_n(x)=0\) and \(g_n(y)=\min\{n,d_X(x,y)\}\) is non-expansive and continuous. By \cref{t:reg-inj}, it extends to a morphism \(f_n\colon X\to[0,n]\). Writing \(P=\prod_{n\in\N}[0,n]\), the product morphism \(f=(f_n)_{n\in\N}\colon X\to P\) satisfies
  \[
    d_P(f(x),f(y))=\sup_{n\in\N}\min\{n,d_X(x,y)\}=d_X(x,y).
  \]
  Therefore, \(P\) is a regular cogenerator by \cref{r:reg-cogen-criterion}.
\end{proof}

      \begin{remark}
        By \cref{r:reg-cogen-criterion} and the proof of \cref{t:reg-cogen}, the interval \([0,\infty]\) is also a regular cogenerator in \(\MetCHsep\) since, for every \(n\in\N\), \([0,n]\) is a subspace of \([0,\infty]\). Similarly, \([0,\infty]^{\circ}\) is a regular cogenerator in \(\MetCHsep\) since, for every \(n\in\N\), \([0,n]^{\circ}\) is a subspace of \([0,\infty]^{\circ}\) and \([0,n]^{\circ}\cong[0,n]\). This observation answers positively the case of the product t-norm of the question posed in \cite{HN18} before Def.~9.3.\footnote{Note that \cite{HN18} considers categories enriched in the quantale \([0,1]\) with multiplication instead of the isomorphic quantale \([0,\infty]\) with addition.} In the terminology of \cite{HN23}, this shows that every separated metric compact Hausdorff space is Priestley (in the metric sense, see \cite[Def.~2.18]{HN23}), answering the question after Corollary~2.22 of \cite{HN23}.
\end{remark}

\begin{theorem}
  The dual of $\MetCHsep$ is equivalent to a variety of (possibly infinitary) algebras. In fact, the representable functor
  \[
    \hom_{\MetCHsep}\biggl(-,\prod_{n \in \N}[0,n]\biggr) \colon \MetCHsep^\op
    \longrightarrow \Set
  \]
  is monadic.
\end{theorem}

\begin{proof}
  By \cite{AH25}, we know already that the dual of $\MetCHsep$ is complete, cocomplete and Barr-exact.
  By \cref{t:reg-inj}, for each $n \in \N$, the object $[0,n]$ is regular injective in $\MetCHsep$.
  Since a product of regular injective objects is regular injective,
  \begin{math}
    \prod_{n \in \N}[0,n]
  \end{math}
  is regular injective, as well. By \cref{t:reg-cogen}, it is also a regular cogenerator in $\MetCHsep$. Therefore, by \cref{d:thm:1}, $\MetCHsepop$ is equivalent to a variety of (possibly infinitary) algebras.
\end{proof}

\section{About the rank of the variety}

We are going to show that $\MetCHsep$ is dually equivalent to a variety of algebras with operations of arity at most countable, and that it is not dually equivalent to any variety of \emph{finitary} algebras; the latter fact was already observed in \cite{AH25}. We strengthen it by showing that $\MetCHsepop$ is not equivalent to any $\mathrm{SP}$-class of finitary algebras, and that it admits no faithful functor to $\Set$ preserving directed colimits.

\begin{theorem} \label{t:char-copres}
  Let \(X\in\MetCHsep\).
  \begin{enumerate}
    \item \label{i:fin-copr}
    \(X\) is finitely copresentable if and only if it has at most one element.

    \item \label{i:aleph1-copr}
    \(X\) is \(\aleph_1\)-copresentable if and only if its underlying topology is metrisable.
  \end{enumerate}
\end{theorem}

\begin{proof}
  \eqref{i:fin-copr}. We first prove the contrapositive of the implication ($\Rightarrow$). Suppose that \(X\) contains two distinct points. Since \(d_X\) is separated, after possibly interchanging them we may choose \(x,y\in X\) such that
  \[
    \alpha\coloneqq d_X(x,y)>0.
  \]
  Choose an increasing sequence \((\beta_n)_{n\geq 1}\) in \([0,\alpha[\) with \(\sup_n\beta_n=\alpha\). Let \(Z_n=\{0,1\}\), with the discrete topology and with the metric given by
  \[
    d_{Z_n}(0,1)=\beta_n,
    \qquad
    d_{Z_n}(1,0)=\infty.
  \]
  For \(m\geq n\), the identity map \(Z_m\to Z_n\) is non-expansive. The limit of the resulting codirected sequence is the discrete two-point space \(Z=\{0,1\}\) with
  \[
    d_Z(0,1)=\alpha,
    \qquad
    d_Z(1,0)=\infty,
  \]
  and all limit projections \(p_n\colon Z\to Z_n\) are the identity on $\{0,1\}$.

  The map \(f\colon Z\to X\), defined by \(f(0)=x\) and \(f(1)=y\), is a morphism. If \(X\) were finitely copresentable, \(f\) would factor through some \(p_n\), say \(f=g\circ p_n\) for a morphism \(g\colon Z_n\to X\). Non-expansiveness of \(g\) would then give
  \[
    \alpha=d_X(x,y)\leq d_{Z_n}(0,1)=\beta_n<\alpha,
  \]
  a contradiction. This proves the implication ($\Rightarrow$).

  Let us now prove the implication ($\Leftarrow$). The one-point space is finitely copresentable because it is terminal. The empty space is finitely copresentable as well: a codirected limit of nonempty compact Hausdorff spaces is nonempty by compactness, and limits in \(\MetCHsep\) are computed as in the category of compact Hausdorff spaces. Hence a codirected limit in \(\MetCHsep\) is empty if and only if one of its components is empty, which shows that \(\hom_{\MetCHsep}(-,\varnothing) \colon \MetCHsepop \to \Set\) preserves directed colimits.

  \eqref{i:aleph1-copr}. This follows from \cite[Corollary~4.38]{HN23}.
\end{proof}

\begin{corollary}
  For each \(n\in\N\), the interval \([0,n]\) is \(\aleph_1\)-copresentable in \(\MetCHsep\). Therefore also
  \begin{math}
    \prod_{n \in \N}[0,n]
  \end{math}
  is so. Hence, \(\MetCHsepop\) is an \(\aleph_1\)-ary variety, that is, it can be presented as a variety with operations of at most countable arity.
\end{corollary}

\begin{corollary}
  \label{cor:not-finitely-accessible}
  The category $\MetCHsepop$ is not finitely accessible.
  Therefore, $\MetCHsepop$ is not equivalent to any variety of finitary algebras (even with an arbitrary set of sorts).
\end{corollary}

\begin{proof}
  By \cref{t:char-copres}, every finitely copresentable object of $\MetCHsep$ has at most one point.
  Since limits in $\MetCHsep$ are computed as in $\Set$, a codirected limit of such objects also has at most one point.
  Thus a two-point space is not a codirected limit of finitely copresentable objects, so its corresponding object in $\MetCHsepop$ is not a directed colimit of finitely presentable objects.
\end{proof}

The following is the analogue of \cite[Cor.~22]{AR20} for separated metric compact Hausdorff spaces.

\begin{corollary}
  The category $\MetCHsepop$ is not equivalent to any $\mathrm{SP}$-class of finitary algebras, that is, any class of finitary algebras closed under subalgebras and arbitrary products, with all homomorphisms between its members.
\end{corollary}

\begin{proof}
  By \cite[Thm.~5.16]{AH25}, $\MetCHsepop$ is Barr-exact, so every internal equivalence relation in it is effective.
  By \cite{Ban83}, an $\mathrm{SP}$-class of finitary algebras with this property is a finitary variety.
  An equivalence with such a class would therefore contradict \cref{cor:not-finitely-accessible}.
\end{proof}

A category is called \emph{finitely concrete} if it admits a faithful functor to \(\Set\) which preserves directed colimits \cite{LRV23}. By an \emph{elementary class of structures} we mean the class of models of a first-order theory in a finitary language; when it is viewed as a category, the morphisms are the homomorphisms.

\begin{proposition}
  \label{prop:not-finitely-concrete}
  The category \(\MetCHsepop\) is not finitely concrete.
  In particular, it is not equivalent to any elementary class of structures.
\end{proposition}

\begin{proof}
  Let $\CH$ be the category of compact Hausdorff spaces and continuous maps between them.
  Let
  \[
    D_\infty\colon\CH\longrightarrow\MetCHsep
  \]
  be the functor that, on objects, equips a compact Hausdorff space with the extended discrete metric, equal to \(0\) on the diagonal and to \(\infty\) off it; on morphisms, it maps a continuous map to itself. Every continuous map between such spaces is non-expansive, so \(D_\infty\) is fully faithful. It preserves products and equalisers: the supremum metric on a product of extended-discrete spaces is again extended discrete, as is the restricted metric on an equaliser. Hence \(D_\infty\) preserves all limits, and \(D_\infty^{\op}\) is a faithful colimit-preserving functor from \(\CH^{\op}\) into \(\MetCHsepop\).

  By Gelfand duality and \cite[proof of Thm.~22]{LRV23}, the category \(\CH^{\op}\) is not finitely concrete (see \cite{Abb26} for a short and direct proof that \(\CH^{\op}\) is not finitely concrete).
  If \(\MetCHsepop\) admitted a faithful directed-colimit-preserving functor to \(\Set\), its composite with \(D_\infty^{\op}\) would contradict this fact.

  For an elementary class of structures with all homomorphisms, the underlying-set functor is faithful and preserves existing directed colimits (see, e.g., \cite[p.~4, after Rem.~5]{LRV23}).
  An equivalence with such a class would therefore give a faithful functor from \(\MetCHsepop\) to \(\Set\) preserving directed colimits.
\end{proof}

\section{The symmetric case}
\label{sec:symmetric-case}

A metric space \(X\) is called \emph{symmetric} whenever the metric \(d\) of \(X\) is symmetric, that is, if
\begin{displaymath}
  d(x,y)=d(y,x)
\end{displaymath}
for all \(x,y\in X\).
We consider now the full subcategory \(\MetCHsepsym\) of \(\MetCHsep\) defined by all symmetric separated metric compact Hausdorff spaces.
We observed in \cite[Thm.~6.1]{AH25} that \(\MetCHsepsym\) is Barr-coexact.
We note that the inclusion functor
\begin{displaymath}
  \MetCHsepsym \longrightarrow \MetCHsep
\end{displaymath}
has a left and a right adjoint (see \cite[Rem.~2.8]{AH25}); the right adjoint \(G \colon \MetCHsep\to\MetCHsepsym\) sends a separated metric compact Hausdorff space to its symmetrisation.
Therefore \(\MetCHsepsym\) is complete and cocomplete, a morphism \(m \colon X\to Y\) in \(\MetCHsepsym\) is a regular monomorphism in \(\MetCHsepsym\) if and only if it is so in \(\MetCHsep\), and \(G(\prod_{n \in \N}[0,n])\) is an \(\aleph_1\)-copresentable regular injective regular cogenerator in \(\MetCHsepsym\).
The same argument as in the proof of \cref{prop:not-finitely-concrete} shows that \(\MetCHsepsymop\) is not finitely concrete.

\begin{corollary}
  \label{cor:symmetric}
  The category \(\MetCHsepsym^{\op}\) is equivalent to a single-sorted \(\aleph _1\)-ary variety. Moreover, it is not finitely concrete and it is not equivalent to any finitary variety, even with an arbitrary set of sorts.
\end{corollary}

\makeatletter
\begingroup
\let\addcontentsline\@gobblethree
\section*{Acknowledgements}
Marco Abbadini is grateful to Enrico Pasqualetto for a discussion thanks to which Marco became aware of Matoušková's paper.
He also thanks Enrico Pasqualetto for explaining how extended metric-topological spaces arose in metric analysis and optimal transport.

\subsection*{Declaration of AI usage}
We used ChatGPT to review parts of the manuscript and to suggest changes to the exposition. Furthermore, ChatGPT called our attention to a small gap in \cite{Mat00} and suggested a correction which we present in \cref{d:lem:1}.

\endgroup
\makeatother

\makeatletter
\begingroup
\let\addcontentsline\@gobblethree
\bibliography{Biblio} 
\bibliographystyle{alphanum2017}
\endgroup
\makeatother

\end{document}